\documentclass[english,a4paper,reqno]{amsart}
\usepackage[mathscr]{eucal}
\usepackage{amsthm,amssymb,amsmath,enumerate}
\usepackage{amscd}
\usepackage{babel}
\usepackage{hyperref}

\newcommand{\dOmega}{{\partial\Omega}}
\newcommand{\veps}{\varepsilon}
\newcommand{\calE}{\mathcal{E}}
\newcommand{\calF}{\mathcal{F}}
\newcommand{\calM}{\mathcal{M}}

\theoremstyle{plain}
\newtheorem{theorem}{Theorem}[section]
\newtheorem{lemma}[theorem]{Lemma}
\newtheorem{proposition}[theorem]{Proposition} 
\newtheorem{corollary}[theorem]{Corollary} 
\theoremstyle{definition}
\newtheorem{definition}[theorem]{Definition} 
\newtheorem{example}[theorem]{Example}
\theoremstyle{remark}
\newtheorem{remark}{Remark}

\numberwithin{equation}{section}

\begin{document}

\title[Monge-Amp\`ere boundary measures]
{Monge-Amp\`ere boundary measures}

\author{Urban Cegrell}
\address{Department of Mathematics and Mathematical Statistics\\
Ume{\aa} University\\
SE-901 87  Ume{\aa}\\
Sweden}
\email{urban.cegrell@math.umu.se}

\author{Berit Kemppe}
\address{Department of Mathematics and Mathematical Statistics\\
Ume{\aa} University\\
SE-901 87  Ume{\aa}\\
Sweden}
\email{berit.kemppe@math.umu.se}





\begin{abstract}
        We study swept-out Monge-Amp\`ere measures
        of plurisubharmonic functions and boundary values 
        related to these measures.
\end{abstract}

\maketitle


\section{Introduction}

\noindent
The purpose of this paper is to study certain boundary measures related to
plurisubharmonic functions on hyperconvex domains.
These measures are obtained as swept-out Monge-Amp\`ere measures
and generalize the boundary measures studied by Demailly in
\cite{Dem1}, see Section \ref{sec_constr}.
A number of properties of the measures, such as
density, support and convergence,
are given in Section \ref{sec_facts}.
The idea is then to use these measures to define and study boundary values of
plurisubharmonic functions on the given domain. This is done in Section
\ref{sec_bdryval}, where we also describe some situations where this coincides
with other notions of boundary values.
Finally in Section \ref{sec_more} we study more general boundary measures
on a more restricted class of hyperconvex domains.
Here we start with a measure on the boundary and find a sequence of
Monge-Amp\`ere measures approximating the given measure.

It is a great pleasure for us to thank Ph\d{a}m Ho\`{a}ng Hi\d{\^{e}}p for
many fruitful comments.


\bigskip\bigskip
\section{Preliminaries}\label{sec_prel}

\noindent
We first recall some definitions needed in this paper.
Let $\Omega$ be a domain in $\mathbb C^n$, $n\ge 2$.
Denote by $PSH(\Omega)$ the plurisubharmonic functions on $\Omega$
and by $PSH^-(\Omega)$ the subclass of nonpositive
functions.
A set $\Omega\subset\mathbb C^n$ is said to be a hyperconvex domain if
it is open, connected and if there exists a function $\varphi\in
PSH^-(\Omega)$ such that
$\{z\in\Omega: \varphi(z)<-c\}\subset\subset\Omega$,
$\forall\, c>0$.
If $\Omega$ is a bounded hyperconvex domain, then it can be shown that the
exhaustion function $\varphi$ can be chosen in $C^\infty(\Omega)\cap
C(\bar\Omega)$ and such that $\int_\Omega (dd^c \varphi)^n<+\infty$
(see \cite{Ceg5}).
This implies for example that the classes defined below are nontrivial.
Unless otherwise stated, $\Omega$ will throughout this paper denote a bounded
hyperconvex domain in $\mathbb C^n$. Also, by a measure we mean a positive
regular Borel measure.

Let $\calE_0(\Omega)$, $\calF(\Omega)$,
$\calE(\Omega)$ and $\calF^{a}(\Omega)$
be the subclasses of $PSH^{-}(\Omega)$ defined as in 
\cite{Ceg1} and \cite{Ceg2}, namely as follows:
\begin{itemize}
\item
        $\calE_0(\Omega)$ is the set of functions
        $u\in PSH(\Omega)\cap L^\infty(\Omega)$ such that
        $\int_\Omega (dd^c u)^n < +\infty$
        and $\lim_{z\to\xi} u(z)=0$, $\forall\,\xi\in\dOmega$
        \smallskip
\item
        $\calF(\Omega)$ is the set of functions
        $u\in PSH(\Omega)$ such that
        there is a sequence $\{u_j\}$ in $\calE_0(\Omega)$
        with the properties that  $u_j\searrow u$ and
        $\sup_j \int_\Omega (dd^c u_j)^n < +\infty$
        \smallskip
\item
        $\calE(\Omega)$ is the set of functions
        $u\in PSH(\Omega)$ such that
        for each $\omega\subset\subset\Omega$ there is
        function $u_\omega\in\calF(\Omega)$ with the properties that
        $u_\omega\ge u$ on $\Omega$ and $u_\omega=u$ on $\omega$
        \smallskip
\item
        $\calF^a(\Omega)$ is the set of functions
        $u\in\calF(\Omega)$ such that
        $\int_E(dd^c u)^n = 0$ for each pluripolar set $E\subset\Omega$
\end{itemize}

For the convenience of the reader, we state some of the
results, concerning these classes, that we use most frequently in this paper.
If nothing else is mentioned, proofs can be found in \cite{Ceg2}.

First, observe that
$PSH^-(\Omega)\cap L^\infty_{\text{loc}}(\Omega)$
is contained in $\calE(\Omega)$ and that
$\calE_0(\Omega) \subset\calF^{a}(\Omega)\subset\calF(\Omega)
\subset\calE(\Omega)$.
The following lemma explains why the functions in $\calE_0(\Omega)$ sometimes are
called \textit{test functions}.

\begin{lemma}\label{l_testfcns}
        If $\varphi\in C_0^\infty(\Omega)$, then there are
        $\phi_1,\phi_2\in\calE_0(\Omega)\cap C(\bar\Omega)$
        such that $\varphi=\varphi_1-\varphi_2.$
\end{lemma}

If $u_1,\ldots,u_n\in\calE(\Omega)$, then $dd^c u_1\wedge\ldots\wedge dd^c u_n$
is defined as the limit measure obtained by combining the following two theorems.

\begin{theorem}\label{th_appr_d}
        Suppose that $u\in PSH^-(\Omega)$. Then there is a sequence
        $\{u_j\}\subset\calE_0(\Omega)\cap C(\bar\Omega)$ such that
        $u_j\searrow u$ on $\Omega$ and
        $\text{supp}\,(dd^c u_j)^n\subset\subset\Omega$ for each $j$.
\end{theorem}

\begin{theorem}\label{th_conv_d}
        For $k=1,\ldots,n$, let $u_k\in\calE(\Omega)$ and
        $\{g_{kj}\}_{j=1}^\infty\subset\calE_0(\Omega)$
        be such that $g_{kj}\searrow u_k$ as $j\to\infty$.
        Then $dd^c g_{1j}\wedge\ldots\wedge dd^c g_{nj}$ is
        weak*-convergent and the limit measure is independent
        on the sequences $\{g_{kj}\}$.
\end{theorem}

A function $u\in\calE(\Omega)$ is a maximal plurisubharmonic function if and only
if $(dd^cu)^n=0$ (see \cite{Blocki} and \cite{Ceg1_2}).
If $u\in\calF(\Omega)$ and $(dd^cu)^n=0$, then $u=0$
(see Theorem 5.15 in \cite{Ceg2}).
Theorem \ref{th_conv_d} can be generalized as follows, see e.g.\ 
Lemma 3.2 in \cite{Ceg00}.

\begin{lemma}\label{l_gen_conv}
        For $k=1,\ldots,n$, let $u_k\in\calE(\Omega)$ and
        $\{g_{kj}\}_{j=1}^\infty\subset\calE(\Omega)$
        be such that $g_{kj}\ge u_k$ and
        $g_{kj}$ tends weakly to $u_k$ as $j\to\infty$.
        If $h\in PSH^-(\Omega)\cap L^\infty(\Omega)$,
        then $h\,dd^c g_{1j}\wedge\ldots\wedge dd^c g_{nj}$ tends weak*
        to $h\,dd^c u_1\wedge\ldots\wedge dd^c u_n$. Moreover,
        if $u_k\in\calF(\Omega)$ then
        $\lim_{j\to\infty}\int_\Omega
        h\,dd^c g_{1j}\wedge\ldots\wedge dd^c g_{nj}=
        \lim_{j\to\infty}\int_\Omega
        h\,dd^c u_1\wedge\ldots\wedge dd^c u_n$.
\end{lemma}

The next lemma contains some useful basic properties of the classes we use.

\begin{lemma}\label{lemma_propert}
        Let $\mathcal{K}\in\{\calE_0,\calF^{a},\calF,\calE\}$,
        then the following holds.
        \begin{enumerate}[(i)]
        \item
                If $u,v\in\mathcal{K}(\Omega)$ and $\alpha,\beta \ge 0$,
                then $\alpha u + \beta v \in\mathcal{K}(\Omega)$.
        \item
                If $u\in\mathcal{K}(\Omega)$ and $v\in PSH^-(\Omega)$,
                then $\text{max}\,\{u,v\}\in\mathcal{K}(\Omega)$.
                In particular, if $u\in\mathcal{K}(\Omega)$,
                $v\in PSH^-(\Omega)$ and $v\ge u$, then
                $v\in\mathcal{K}(\Omega)$.
        \end{enumerate}
\end{lemma}

Note that functions in $\calF(\Omega)$ have finite total Monge-Amp\`ere
mass.
Also, they have in some sense boundary values zero, which can be seen e.g.\ in
the following formula for partial integration.

\begin{theorem}\label{th_pi}
        Let $v,u_1,\ldots,u_n\in\calF(\Omega)$. Then
        \begin{equation*}
                \int_\Omega v\,dd^c u_1\wedge dd^c u_2\wedge\ldots
                \wedge dd^c u_n = 
                \int_\Omega u_1\,dd^c v\wedge dd^c u_2\wedge\ldots
                \wedge dd^c u_n.
        \end{equation*}
\end{theorem}

Since bounded function cannot put Monge-Amp\`ere mass on pluripolar sets
(see e.g.\ \cite{B&T}), we
have that $\calF(\Omega)\cap L^\infty(\Omega)\subset\calF^{a}(\Omega)$.
Moreover, Theorem 5.5 and Theorem 5.8 in \cite{Ceg2} gives:

\begin{lemma}\label{l_pp}
        If $u_1,\ldots,u_{n-1}\in\calF(\Omega)$ and $v\in\calF^{a}(\Omega)$
        or $v\in PSH^-(\Omega)\cap L^\infty(\Omega)$,
        then $dd^c u_1\wedge\ldots \wedge dd^c u_{n-1}\wedge dd^c v$
        vanishes on pluripolar sets.
\end{lemma}

We conclude this section with some notation needed in this paper.
Let $\Omega$ and $u\in\calE(\Omega)$ be given and choose a fundamental sequence
$\{\Omega_j\}$ of strictly pseudoconvex domains,
i.e.\
$\Omega_j\subset\subset\Omega_{j+1}\subset\subset\Omega$ and
$\cup_{j=1}^\infty\Omega_j=\Omega$.
For each $j$ define
\begin{equation}\label{uj_def}
        u^j = \sup\,\{\varphi\in PSH(\Omega):
        \varphi|_{\Omega\setminus\Omega_j}\le
        u|_{\Omega\setminus\Omega_j}\}.
\end{equation}
Note that since $\Omega_j$ has $C^2$ boundary, it follows that
$u^j=(u^j)^*$, the smallest upper
semicontinuous majorant of $u^j$, so $u^j$ is plurisubharmonic.
Moreover $u\le u^j \le u^{j+1}\le 0$, so each $u^j\in\calE(\Omega)$
and the same holds for $\tilde u = (\lim_{j\to\infty} u^j)^{*}$.
It follows that
$\tilde u$ is the smallest maximal plurisubharmonic majorant of $u$
and that $\tilde u$ is independent of the chosen sequence $\{\Omega_j\}$.
In \cite{Ceg00} the following classes were defined:
\begin{equation*}\begin{split}
        \mathcal{N}(\Omega) & =
        \{u\in\calE(\Omega): \tilde u = 0\}\\
        \calM(\Omega) & = 
        \{u\in \calE(\Omega): (dd^c u)^n = 0\}
\end{split}\end{equation*}
Thus $\calM(\Omega)$ is the class of maximal
plurisubharmonic functions in $\calE(\Omega)$.
Note that $\mathcal{N}(\Omega)$ contains $\calF(\Omega)$, since if
$u\in\calF(\Omega)$, then $\tilde u$ is a maximal function in $\calF(\Omega)$
so $\tilde u=0$.
It also follows that if $u\in\calF(\Omega)$, then $u^j\nearrow
0$ outside a pluripolar subset of $\Omega$ (see \cite{Lel2} or \cite{B&T}).

Finally, we say that $u\in\calE(\Omega)$ has boundary values $\tilde u$ if
there is a function $\psi\in\mathcal{N}(\Omega)$ such that
$\tilde u\ge u \ge \tilde u + \psi$.
Given $H\in\calM(\Omega)$ we define
$$\calF(\Omega,H) =
\{u\in PSH(\Omega): H\ge u \ge H+\psi,\,\psi\in \calF(\Omega)\},$$
which is a subclass of $\calE(\Omega)$.
It follows that if $u\in\calF(\Omega,H)$ then
$\tilde u = H$. Also, $\calF(\Omega,0)=\calF(\Omega)$.



\bigskip\bigskip
\section{Construction of the boundary measures $\mu_u$}\label{sec_constr}

\noindent
In this section we show that every function in $\calF(\Omega)$ gives rise
to a measure on the boundary of $\Omega.$
Let $u\in\calF(\Omega)$ be given, choose a fundamental sequence
$\{\Omega_j\}$ of strictly pseudoconvex domains
and let $u^j$ be defined by (\ref{uj_def}).
Then $u\le u^j \le u^{j+1}\le 0$, so each $u^j\in\calF(\Omega)$.
Moreover, Stokes' theorem implies that
$\int_\Omega(dd^c u^j)^n = \int_\Omega(dd^c u)^n<+\infty$, and by maximality
$(dd^c u^j)^n$ is concentrated on $\Omega\setminus\Omega_j$.

\begin{theorem}\label{th1}
        Suppose that $u\in\calF(\Omega)$.
        Then $\{(dd^cu^j)^n\}$ is a weak*-convergent
        sequence, which defines a positive measure
        $\mu_u$ on $\dOmega$.
        Also $\lim_{j\to\infty}\int_\Omega\varphi\,(dd^cu^j)^n$ exists for all
        $\varphi\in PSH(\Omega)\cap L^{\infty}(\Omega)$.
\end{theorem}

\begin{proof}
        Choose W to be a strictly pseudoconvex set containing the closure of 
        $\Omega$.
        First assume that $\varphi\in PSH(\Omega)\cap L^\infty(\Omega)$
        and $\varphi\le 0$, then
        \begin{equation}\label{ineq_1}
                -\infty<\int_\Omega \varphi\,(dd^c u)^n
                \le \int_\Omega \varphi\,(dd^c u^j)^n
                \le \int_\Omega \varphi\,(dd^c u^{j+1})^n
                \le \sup_\Omega \varphi\int_\Omega (dd^c u)^n.
        \end{equation}
        To see this, approximate $\varphi$ with functions in
        $\calE_0(\Omega)$ and use
        partial integration in $\calF(\Omega)$ (see Section \ref{sec_prel}).
        Since all Monge-Amp\`ere measures involved have the same total mass,
        it follows that (\ref{ineq_1}) holds
        for all $\varphi\in PSH(\Omega)\cap L^\infty(\Omega)$.
        Thus $\{\int_\Omega \varphi\,(dd^c u^j)^n\}$ is a bounded monotone
        sequence, so $\lim_{j\to\infty}\int_\Omega \varphi\,(dd^c u^j)^n$ exists for all
        $\varphi\in PSH(\Omega)\cap L^\infty(\Omega)$.
        In particular the limit exists for
        $\varphi\in C_0^\infty(W)$
        (see Lemma \ref{l_testfcns}).
        Since each $(dd^c u^j)^n$ is a positive distribution
        on $C_0^\infty (W)$, it follows from standard distribution theory that
        the convergence
        in fact holds for all $\varphi\in C_0(W)$. Also
        the limit distribution itself is
        positive and thus defines
        a positive regular Borel measure $\mu_u$ on $W$, which by the construction is
        concentrated on $\dOmega$.
\end{proof}

In this manner we may, to each
$u\in\calF(\Omega)$, associate a positive measure $\mu_u$, and it follows for
example that
\begin{equation}\label{limit_1}
        \int_\dOmega \varphi\, d\mu_u =
        \lim_{j\to\infty}\int_\Omega \varphi\,(dd^c u^j)^n
\end{equation}
holds
for all $\varphi\in C_0(W)$, in particular for
$\varphi\in C(\bar\Omega)$.
We also have that $\int_\dOmega d\mu = \int_\Omega (dd^c u)^n$, which implies
that $\mu_u=0$ if and only if $u=0$ (since $u\in\calF(\Omega)$).
Note that $\mu_u$ does not depend on the chosen sequence $\{\Omega_j\}$.
Note also that by applying (\ref{ineq_1}) to
 $\varphi$ and $-\varphi$ we get that
\begin{equation}\label{eq_1}
        \int_\Omega \varphi\,(dd^c u^j)^n=
        \int_\Omega \varphi\,(dd^c u)^n,\quad
        \forall\,\varphi\in PH(\Omega)\cap L^\infty(\Omega),
\end{equation}
where $PH(\Omega)$ denotes the pluriharmonic functions on $\Omega$.


In \cite{Dem1} Demailly defines a set
of Monge-Amp\`ere boundary measures in the following setting.
Let $X$ be a Stein manifold of dimension $n$ and $\Omega\subset\subset X$ an
open hyperconvex subset.
Assume that $\phi:\Omega\to [-\infty,0)$ is a continuous plurisubharmonic
exhaustion function such that $\int_\Omega (dd^c \phi)^n < +\infty$.
For each $r<0$ define:
\begin{eqnarray*}
        B(r)&=&\{z\in\Omega : \phi(z) < r\}\\
        S(r)&=&\{z\in\Omega : \phi(z) = r\}\\
        \phi_r(z)&=&\max\,\{\phi(z),r\}
\end{eqnarray*}
It is then shown that
\begin{equation}
         (dd^c \phi_r)^n =
         \chi_{\Omega\setminus B(r)}\cdot  (dd^c \phi)^n +
         \mu_{\phi,r} 
\end{equation}
where $\mu_{\phi,r}$ is a positive measure concentrated on $S(r)$.
Furthermore, when $r\to 0$ then $\mu_{\phi,r}$ converges in a weak sense
to a positive measure $\tilde{\mu}_\phi$ concentrated on $\dOmega$.
(More explicitly it is shown that
$\lim_{r\to 0}\int h\,d\mu_{\phi,r}$ exists
$\forall\, h\in C^2(X,\mathbb{R})$.)

Now consider the case when $X=\mathbb{C}^n$, then the function
$\phi$ is in $\calF(\Omega)$ so we can define $\mu_\phi$ according to
Theorem \ref{th1}. Choose a
sequence $\{r_j\}$ such that
$r_j\nearrow 0$ and let $\Omega_j = B(r_j)$. Then $ \phi_{r_j}=\max\,\{\phi,r_j\}$
is equal to the function $\phi^j$ defined as in (\ref{uj_def}).
Note that $\Omega_j$ is not necessarily
strictly pseudoconvex in this setting, only hyperconvex.
However, this is enough in the proof of
Theorem \ref{th1}, since we only use the smoothness of $\dOmega_j$ to ensure that
the function $\phi^j$ is plurisubharmonic.
Hence
\begin{equation}
        (dd^c \phi_{r_j})^n =
         \chi_{\Omega\setminus B(r_j)}\cdot  (dd^c \phi)^n +
         \mu_{\phi,r_j},
\end{equation}
where the left hand side converges to the boundary measure $\mu_\phi$
and the right hand side to
$0+\tilde{\mu}_\phi$ (since $\int_\Omega (dd^c \phi)^n < +\infty$).
This shows that $\mu_\phi = \tilde{\mu}_\phi$, so in particular Demailly's
boundary measures
form a subset of those defined in Theorem \ref{th1}, when $X=\mathbb{C}^n$.

Also, note that if $u\in\calE_0(\Omega)\cap C(\bar\Omega)$ then $u$
satisfies
the conditions in Demailly's definition, so for boundary measures corresponding
to such functions we may use Demailly's results.

The following theorem, where $u^j$ is defined by (\ref{uj_def}), generalizes a
formula considered by Demailly in \cite{Dem1}.

\begin{theorem}\label{th_formula}
        Assume that $u\in \calF(\Omega)$, $h\in \calE(\Omega)$,
        $\int_\Omega h\,(dd^c u)^n > -\infty $ and
        that $dd^ch\wedge (dd^c u)^{n-1}$ vanishes on pluripolar sets. 
        Then 
        $$\lim_{j\to\infty} \int_\Omega h\,(dd^c u^j)^n =
        \int_\Omega h\,(dd^c u)^n - \int_\Omega u\,dd^c h\wedge (dd^c u)^{n-1}.$$
\end{theorem}

Note that the conditions in this theorem are satisfied if for example
$u\in\calF(\Omega)$ and $h\in PSH^{-}(\Omega)\cap L^\infty(\Omega)$
(see Lemma \ref{l_pp}).
Actually, it is enough that $h\in PSH(\Omega)\cap L^\infty(\Omega)$,
since $\int_\Omega(dd^c u^j)^n = \int_\Omega(dd^c u)^n$.

\begin{proof}[Proof of Theorem \ref{th_formula}]
        First we claim the following.
        \begin{enumerate}[(i)]
        \item\label{c1}
                $\displaystyle
                \int_\Omega u\,dd^c h\wedge(dd^c u)^{n-1} > -\infty$
        \item\label{c2}
                $\displaystyle\lim_{j\to\infty}
                \int_\Omega u^j\,dd^c h\wedge(dd^c u)^{n-1}=0$
        \item\label{c3}
                $\displaystyle
                \int_\Omega h\,(dd^c u^j)^{n-p+1}\wedge (dd^c u)^{p-1}\ge
                \int_\Omega h\,(dd^c u^j)^{n-p}\wedge (dd^c u)^{p}\ge
                \int_\Omega h\,(dd^c u)^{n}$,
                $1\le p\le n-1$
        \item\label{c4}
                $\displaystyle
                \int_\Omega
                h\,dd^c(u^j-u)\wedge(dd^c u^j)^{n-p}\wedge(dd^c u)^{p-1}=$
                
                \noindent
                $\displaystyle =\int_\Omega u^j\,dd^c h\wedge dd^c(u^j-u)
                \wedge(dd^c u^j)^{n-p-1}\wedge(dd^c u)^{p-1}=$
                
                \noindent
                $\displaystyle =\int_\Omega
                (u^j-u)\,dd^c h\wedge(dd^c u^j)^{n-p}\wedge(dd^c u)^{p-1}\ge0$,
                $1\le p\le n$
        \end{enumerate}
        For the proof of (\ref{c1}), choose a sequence $\{h_k\}$ in
        $\calE_0(\Omega)$ decreasing to $h$ on $\Omega$.
        Then $dd^c h_k\wedge(dd^c u)^{n-1}$ converges weak* to
        $dd^c h\wedge(dd^c u)^{n-1}$ (Lemma \ref{l_gen_conv}).
        Combining this with the fact that $u$ is upper semicontinuous it
        follows that
        \begin{equation*}
        \begin{split}
                \int_\Omega (-u)\,dd^c h\wedge(dd^c u)^{n-1}&\le
                \limsup_{k\to\infty}
                \int_\Omega (-u)\,dd^c h_k\wedge(dd^c u)^{n-1}=\\
                & =\limsup_{k\to\infty}\int_\Omega(-h_k)\,(dd^c u)^{n}
                 = \int_\Omega (-h)\,(dd^c u)^{n}<+\infty
        \end{split}
        \end{equation*}
        (where we have used partial integration in $\calF(\Omega)$).
        Since $u^j\nearrow 0$ outside a pluripolar set
        (see Section \ref{sec_prel})
        and since $dd^ch\wedge (dd^c u)^{n-1}$ puts no mass there,
        (\ref{c1}) implies (\ref{c2}) by dominated convergence.
        To see (\ref{c3}), use the same technique as in Theorem \ref{th1}.
        Finally (\ref{c4}) follows from partial integration, using the fact
        that $h$ is locally in $\calF(\Omega)$ and that $u^j-u$ is compactly
        supported in $\Omega$. This proves the claim.
        
        Now using (\ref{c4}) we have that
        \begin{equation*}
        \begin{split}
                \int_\Omega u\,dd^c h\wedge(dd^c u)^{n-1} &=
                \int_\Omega (u-u^j)\,dd^c h\wedge(dd^c u)^{n-1} +
                \int_\Omega u^j\,dd^c h\wedge(dd^c u)^{n-1}=\\
                & = \int_\Omega h\,dd^c(u-u^j)\wedge(dd^c u)^{n-1} +
                \int_\Omega u^j\,dd^c h\wedge(dd^c u)^{n-1},
        \end{split}
        \end{equation*}
        so we can write
        \begin{multline*}
                \int_\Omega h\,(dd^c u^j)^n - \int_\Omega h\,(dd^c u)^n +
                \int_\Omega u\,dd^c h\wedge (dd^c u)^{n-1} =\\
                = \int_\Omega h\,(dd^c u^j)^n -
                \int_\Omega h\,dd^c u^j\wedge(dd^c u)^{n-1} +
                \int_\Omega u^j\,dd^c h\wedge(dd^c u)^{n-1}
        \end{multline*}
        where the last integral tends to $0$ according to (\ref{c2}).
        Moreover
        \begin{multline*}
                \int_\Omega h\,(dd^c u^j)^n -
                \int_\Omega h\,dd^c u^j\wedge(dd^c u)^{n-1} =\\
                =\sum_{p=1}^{n-1}
                \left(\int_\Omega h\,(dd^c u^j)^{n-p+1}\wedge(dd^c u)^{p-1}-
                \int_\Omega h\,(dd^c u^j)^{n-p}\wedge(dd^c u)^{p}\right) =
                \sum_{p=1}^{n-1} a_p
        \end{multline*}
        where each $a_p\ge 0$ by (\ref{c3}). Using (\ref{c4}) we have that
        \begin{equation*}
        \begin{split}
                a_p & = \int_\Omega
                h\,dd^c(u^j-u)\wedge(dd^c u^j)^{n-p}\wedge(dd^c u)^{p-1}=\\
                & = \int_\Omega u^j\,dd^c h\wedge dd^c(u^j-u)
                \wedge(dd^c u^j)^{n-p-1}\wedge(dd^c u)^{p-1}\le\\
                & \le -\int_\Omega
                u^j\,dd^c h\wedge(dd^c u^j)^{n-p-1}\wedge(dd^c u)^p.
        \end{split}
        \end{equation*}
        Now, the second expression in (\ref{c4}) implies that
        $\int_\Omega u^j\,dd^c h\wedge(dd^c u^j)^{n-k}\wedge(dd^c u)^{k-1}$
        is decreasing in $k$, so it follows that
        $0\le a_p\le -\int_\Omega u^j\,dd^c h\wedge(dd^c u)^{n-1}$.
        Hence (\ref{c2}) implies that each term $a_p\to 0$ as $j\to\infty$
        and the theorem is proved.
 \end{proof}
 
 \begin{remark}\label{rem_pi}
Combining the preceeding theorem with (\ref{limit_1}), we have the following
formula. Given $u\in\calF(\Omega)$,
\begin{equation}\label{pi_formula}
        \int_\Omega h\,(dd^c u)^n=
        \int_\Omega u\,dd^c h\wedge (dd^c u)^{n-1} +
        \int_\dOmega h\,d\mu_u,
        \quad\forall\,h\in PSH(\Omega)\cap C(\bar\Omega).
\end{equation}
In Section \ref{sec_facts} (Corollary \ref{cor_supp_S}) we will show that there
is a set $S\subset\dOmega$
such that $\text{supp}\,\mu_u=S$ for each $u\in\calF(\Omega)$, $u\ne 0$.
Hence (\ref{pi_formula}) gives a partial integration formula for
$h\in PSH^-(\Omega)\cap C(\bar\Omega)$ such that $h|_S=0$.
From Theorem \ref{th_equal} in Section \ref{sec_bdryval} it follows that if
$u\in\calF^{a}(\Omega)$, then (\ref{pi_formula}) is valid for
$h\in PSH(W)\cap L^\infty(W)$, where $W$ is some neighbourhood of $\Omega$.

We also get a Jensen-type inequality; given $u\in\calF(\Omega)$,
\begin{equation}\label{J_ineq}
        \int_\Omega h\,(dd^c u)^n \le
        \int_\dOmega h\,d\mu_u,
        \quad\forall\,h\in PSH(\Omega)\cap C(\bar\Omega).
\end{equation}
If $h\in PSH(W)$ for some neighbourhood $W$ of $\Omega$, then using
convolution we may find functions $h_k\in PSH(W')\cap C(W')$, where
$\bar\Omega\subset W'\subset\subset W$, such that $h_k\searrow h$ on $W'$.
Therefore (\ref{J_ineq}) holds true if $h\in PSH(W)$ and $u\in\calF(\Omega)$.
 \end{remark}


\bigskip\bigskip
\section{Some properties of the boundary measures $\mu_u$}\label{sec_facts}

\noindent
In this section we investigate some properties of the boundary
measures $\mu_u$ defined in Section \ref{sec_constr}.
Recall that a hyperconvex domain $\Omega$ is called B-regular if each
continuous function on $\dOmega$ can be extended continuously to a
plurisubharmonic function on $\Omega$ (see \cite{Sib}).

\begin{theorem}\label{th3}
        Let $\mu$ be a finite positive measure on $\partial\Omega$,
        where $\Omega$ is a
        bounded  B-regular domain.
        Then $\mu$ is in the weak* closure of
        $\{\mu_u : u\in \calF(\Omega)\}$.
\end{theorem}

\begin{proof}
        For simplicity, assume that $\mu(\partial\Omega)=1$.
        Choose a sequence of measures
        $$\mu_k=\sum_{j=1}^{N_k}a_j^k\delta_{z_j^k},
       \text{ where }
        {\{z_j^k\}}_{j=1}^{N_k}\subset\Omega
        \text{ and }
        \sum_{j=1}^{N_k}a_j^k=1$$
        such that
        \begin{equation}\label{th3_eq1}
                \lim_{k\to\infty}\int_\Omega h\,d\mu_k =
                \int_\dOmega h\,d\mu,
                \quad\forall\,h\in C(\bar\Omega).
        \end{equation}
        Let e.g.\ $a_j^k=\mu(A_j^k)$ and $z_j^k\in A_j^k\cap\Omega$,
        where ${\{A_j^k\}}_{j=1}^{N_k}$ is a partition of $\bar\Omega$
        such that $\text{diam}({A_j^k})\le\frac{1}{2^k}$,
        and use the fact that $h$ is uniformly continuous on $\bar\Omega$.
        For each $k$, consider $g_k(z)$, the multipole pluricomplex Green's
        function for $\Omega$ with poles at $\{z_j^k\}$ with weights
        $\{(a_j^k)^{1/n}\}$ (see \cite{Lel} and \cite{Pol}).
        Then $g_k\in \calF(\Omega)$ and $(dd^c g_k)^n = \mu_k$.
        Form $\tilde\mu_k = \lim_{i\to\infty}(dd^c(g_k)^i)^n$ as in section
        \ref{sec_constr}.
        Then for each $k$
        \begin{equation}\label{th3_eq2}
                \int_\dOmega d\tilde\mu_k = \int_\Omega(dd^c g_k)^n =
                \int_\Omega d\mu_k = 1 = \int_\dOmega d\mu
        \end{equation}
        and from (\ref{limit_1}) and (\ref{ineq_1}) it follows that
        \begin{equation}\label{th3_eq3}
                \int_\dOmega \varphi\,d\tilde\mu_k =
                \lim_{i\to\infty}\int_\Omega \varphi\,(dd^c (g_k)^i)^n \ge
                 \int_\Omega \varphi\,(dd^c g_k)^n =
                 \int_\Omega \varphi\,d\mu_k
        \end{equation}
        for $\varphi\in PSH(\Omega)\cap C(\bar\Omega)$.
        Let $\{\tilde\mu_{k_m}\}$ be any weak*-convergent subsequence
        of $\{\tilde\mu_k\}$.
        (Such a subsequence exists since the measures $\{\tilde\mu_k\}$
        have uniformly bounded total mass.)
        Now let $t\in C(\dOmega)$, $t\le 0$ be given. 
        Since $\Omega $ is B-regular there is a
        $\varphi\in PSH(\Omega)\cap C(\bar\Omega)$ with $\varphi = t$
        on $\dOmega$.
        Hence, by (\ref{th3_eq1}) and (\ref{th3_eq3}),
        $$
        \int_\dOmega t\,d\mu = \lim_{m\to\infty}\int_\Omega \varphi\,d\mu_{k_m}
        \le \lim_{m\to\infty}\int_\dOmega\varphi\,d\tilde\mu_{k_m} =
        \lim_{m\to\infty}\int_\dOmega t \,d\tilde\mu_{k_m}.
        $$
        This shows that $\mu\ge\lim_{m\to\infty}\tilde\mu_{k_m}$.
        It then follows from (\ref{th3_eq2}) that they have the same total
        mass, so $\mu=\lim_{m\to\infty}\tilde\mu_{k_m}$ and the theorem is proved.
        Note that since the argument is valid for any weak*-convergent
        subsequence,
        it follows that $\{\tilde\mu_k\}$ itself tends
        weak* to $\mu$.
\end{proof}

Later in this section, we will show that not every positive measure on
$\dOmega$ is in $\{\mu_u: u\in\calF(\Omega)\}$, see for example Proposition
\ref{prop_support}.
Moreover, the assumption of B-regularity cannot be removed in
Theorem \ref{th3}, see for example Corollary \ref{cor_supp_S} and
Example \ref{ex_bidisc}.
Before we can prove this, we need the following
convergence property.

\begin{proposition}\label{appr_prop}
        Suppose that $u\in\calF(\Omega)$ and that $\{u_k\}$ is a decreasing
        sequence in $\calF(\Omega)$ such that $u_k\searrow u$ on $\Omega$.
        Then $\mu_{u_k}$ converges weak* to $\mu_u$.
\end{proposition}

\begin{proof}
        Let $h\in\calE_0(\Omega')\cap C(\bar\Omega')$ where
        $\Omega'\supset\bar\Omega$. Then (\ref{pi_formula})
        gives that
        $$\int_\dOmega h\,d\mu_u =
         \int_\Omega h\,(dd^c u)^n -
        \int_\Omega u\,dd^c h\wedge(dd^c u)^{n-1}$$
        and that for each $k$
        $$\int_\dOmega h\,d\mu_{u_k} =
        \int_\Omega h\,(dd^c u_k)^n -
        \int_\Omega u_k\,dd^c h\wedge(dd^c u_k)^{n-1}.$$
        From Lemma \ref{l_gen_conv} it follows that
        $\lim_{k\to\infty} \int_\Omega h\,(dd^c u_k)^n =
         \int_\Omega h\,(dd^c u)^n$.
        Moreover,
        $\lim_{k\to\infty}\int_\Omega u_k\,dd^c h\wedge(dd^c u_k)^{n-1}
        =\int_\Omega u\,dd^c h\wedge(dd^c u)^{n-1}$ by the following
        calculations.
        Since $u\le u_k$ for each $k$, Lemma 3.3 in \cite{ACCP}
        implies that
        $$\int_\Omega u\,dd^c h\wedge(dd^c u)^{n-1}\le
        \int_\Omega u\,dd^c h\wedge(dd^c u_k)^{n-1}\le
        \int_\Omega u_k\,dd^c h\wedge(dd^c u_k)^{n-1}$$
        for each $k$. Hence, for fixed $k_0$,
        \begin{equation*}
        \begin{split}
                \int_\Omega u\,dd^c h\wedge(dd^c u)^{n-1}
                & \le \liminf_{k\to\infty}
                \int_\Omega u_k\,dd^c h\wedge(dd^c u_k)^{n-1}\\
                & \le\limsup_{k\to\infty}
                \int_\Omega u_k\,dd^c h\wedge(dd^c u_k)^{n-1}\\
                & \le\limsup_{k\to\infty}
                \int_\Omega u_{k_0}\,dd^c h\wedge(dd^c u_k)^{n-1}\\
                & \le \int_\Omega u_{k_0}\,dd^c h\wedge(dd^c u)^{n-1},
        \end{split}
        \end{equation*}
        where the last inequality follows since
        $dd^c h\wedge(dd^c u_k)^{n-1}$ is weak*-convergent to
        $dd^c h\wedge(dd^c u)^{n-1}$ (Lemma \ref{l_gen_conv}) and
        $u_{k_0}$ is upper semicontinuous.
        Now, the claim follows if we let $k_0\to\infty$.
        
        Thus
        \begin{equation}\label{limit_h}
                \lim_{k\to\infty}\int_\dOmega h\,d\mu_{u_k} =
                \int_\dOmega h\,d\mu_{u}
        \end{equation}
        holds true for $h\in\calE_0(\Omega')\cap C(\bar\Omega')$
        and therefore for $h\in C_0^\infty(\Omega')$.
        By standard distribution theory it follows that (\ref{limit_h}) holds
        for $h\in C_0(\Omega')$ and hence for $h\in C(\dOmega)$.
\end{proof}

Recall from Section \ref{sec_constr} that for functions in
$\calE_0(\Omega)\cap C(\bar\Omega)$ we can apply the results of Demailly in
\cite{Dem1}. We make use of this fact in the proof of the following
proposition.

\begin{proposition}\label{ineq_prop}
        If $u$ and $v$ are functions in $\calF(\Omega)$
        such that $u\le v$, then $\mu_u\ge\mu_v$.
\end{proposition}

\begin{proof}
        Take $\{u_k\},\{w_k\}\subset
        \calE_0(\Omega)\cap C(\bar\Omega)$ such that
        $u_k\searrow u$ and $w_k\searrow v$.
        Let $v_k=\max\,\{u_k,w_k\}$. Then $v_k\in
        \calE_0(\Omega)\cap C(\bar\Omega)$,
        $v_k\searrow v$ and $u_k\le v_k$.
        By Theorem 3.4 in \cite{Dem1}
        $\mu_{u_k}\ge\mu_{v_k}$ for each $k$.
        Using Proposition \ref{appr_prop} it follows that
        $\mu_u\ge\mu_v$.
\end{proof}

\begin{remark}
        When $\Omega$ is B-regular there is a slightly more direct proof
        of Proposition \ref{ineq_prop},
        not using Demailly's results. If in that case $f\in C(\dOmega)$, $f\le 0$
        is given, it may be extended to a function in
        $PSH^-(\Omega)\cap C(\bar\Omega)$. Since $u\le v$ we have that
        $u^j\le v^j$ for each $j$, which (see the proof of
        Theorem \ref{th1}) implies that
        $\int_\Omega f\,(dd^c u^j)^n\le \int_\Omega f\,(dd^c v^j)^n$
        for each $j$. From  (\ref{limit_1}) it follows that
        $\int_\dOmega f\,d\mu_u\le\int_\dOmega f\,d\mu_v$,
        so we have, by the regularity of $\mu_u$ and $\mu_v$,
        that $\mu_u\ge\mu_v$.
\end{remark}

\begin{corollary}\label{cor_max}
        Suppose that $u\in \calF(\Omega)$, then $\mu_u = \mu_{\max\,\{u,-1\}}$.
\end{corollary}

\begin{proof}
        Let $v = \max\,\{u,-1\}$, then $\mu_u\ge \mu_v$ by Proposition
        \ref{ineq_prop}. Take $\{u_k\}\subset\calE_0(\Omega)$ such that
        $u_k\searrow u$ and let $v_k = \max\,\{u_k,-1\}$.
        Then $v_k\in\calE_0(\Omega)$, $v_k\searrow v$ and
        $v_k = u_k$ on $\Omega\setminus\{u_k<-1\}$
        (note that $\{u_k<-1\}\subset\subset\Omega$).
        Using Theorem 5.1 in \cite{Ceg2} and Stokes theorem,
        it follows that
        \begin{multline*}
                \int_\dOmega d\mu_u =
                \int_\Omega (dd^c u)^n =
                \lim_{k\to\infty}\int_\Omega (dd^c u_k)^n =\\
                =\lim_{k\to\infty}\int_\Omega (dd^c v_k)^n =
                \int_\Omega (dd^c v)^n
                = \int_\dOmega d\mu_v,
        \end{multline*}
         so $\mu_u = \mu_v$.
\end{proof}

We will now use this corollary to show that each $\mu_u$ vanishes on
pluripolar sets. We start with two technical lemmas.

\begin{lemma}\label{tec_lemma}
        Suppose that $u\in\calF(\Omega)$ and that $\varphi$ is in
        $PSH(\Omega)\cap L^\infty(\Omega)$
        and upper semicontinuous on some neighbourhood of $\bar\Omega$.
        Then
        $$\lim_{j\to\infty}\int_\Omega \varphi\,(dd^c u^j)^n \le
        \int_\dOmega\varphi\,d\mu_u.$$
\end{lemma}
 
\begin{proof}[Proof of Lemma \ref{tec_lemma}]
        Choose $\Omega'$ and $\Omega''$ such that $\varphi$ is
        upper semicontinuous on $\Omega'$ and
        $\Omega\subset\subset\Omega''\subset\subset\Omega'$.
        Then there is a decreasing sequence
        $\{\varphi_k\}$ of continuous functions on $\Omega''$ that are
        bounded above and that converge to $\varphi$ on
        $\Omega''$.
        Using equality (\ref{limit_1}) we have
        that
        $$\lim_{j\to\infty}\int_\Omega \varphi\,(dd^c u^j)^n\le
        \lim_{j\to\infty}\int_\Omega \varphi_k\,(dd^c u^j)^n=
        \int_\dOmega\varphi_k\,d\mu_u$$
        for each $k$.
        Hence the lemma follows by letting $k\to\infty$.
\end{proof}

\begin{lemma}\label{lemma_support}
        Let $E\subset\dOmega$ be a pluripolar set
        and $u\in\calF(\Omega)$. Suppose that there is a function
        $g\in PSH(\Omega')$, where $\Omega'\supset\Omega$,
        such that $E\subset\mathcal{S}_g=\{z: g(z)=-\infty\}$
        and $(dd^c u)^n$ is concentrated on $\Omega\setminus\mathcal{S}_g$.
        Then $\mu_u(E)=0$.
\end{lemma}

\begin{proof}
        By subtracting a suitable constant
        we may assume that $g\le 0$ on $\bar\Omega$.
        For each positive integer $k$, define
        $h_k=\max\,\{\frac{1}{k}\cdot g, -1\}$. Then
        from ($\ref{ineq_1}$) and Lemma \ref{tec_lemma} it follows that
        $$-\infty<\int_\Omega h_k\,(dd^c u)^n \le
        \lim_{j\to\infty} \int_\Omega h_k\,(dd^c u^j)^n\le
        \int_\dOmega h_k\,d\mu_u\le
        \int_E h_k\,d\mu_u = -\mu_u(E),$$
        since $h_k\le 0$ on $\bar\Omega$ and $h_k=-1$ on $E$.
        Moreover,
        $h_k (z)\nearrow 0$ for all
        $z\in\Omega\setminus\mathcal{S}_g$, as $k\to\infty$,
        so
        $\lim_{k\to\infty}\int_\Omega h_k\,(dd^c u)^n = 0$.
        Hence $\mu_u(E)=0$.
\end{proof}

\begin{proposition}\label{prop_support}
        If $u\in\calF(\Omega)$, then $\mu_u$ vanishes on
        pluripolar subsets of $\dOmega$.
\end{proposition}

\begin{proof}
        If $u\in\calF(\Omega)$ then $v = \max\,\{u,-1\}\in\calF^{a}(\Omega)$
        and from Corollary \ref{cor_max} we know that $\mu_u=\mu_v$.
        Now, for functions in
        $\calF^{a}(\Omega)$ the conditions in Lemma \ref{lemma_support} are
        satisfied for each pluripolar set $E\subset\dOmega$, so the
        proposition follows.
\end{proof}

The next proposition enables us to say more about the
support of the $\mu_u$-measures.

\begin{proposition}\label{prop_compare}
        Assume that $u,v\in\calE_0(\Omega)$ are
        strictly negative functions such that
        $\text{supp}\,(dd^c u)^n\subset\subset\Omega$ and
        $\text{supp}\,(dd^c v)^n\subset\subset\Omega$.
        Then there are constants $a,b > 0$ such that
        $$a\mu_u\le\mu_v\le b\mu_u.$$
        In particular, $\text{supp}\,\mu_u = \text{supp}\,\mu_v$.
\end{proposition}

\begin{lemma}\label{lemma1}
        Assume that $u\in\calF^{a}(\Omega)$, $u\ne0$, $v\in\mathcal E(\Omega)$
        and that $u \geq v$ on $\text{supp}\,(dd^cu)^n$.
        Then $u \geq v$ on $\Omega$.
\end{lemma}

\begin{proof}
        Assume that $u(z_0)<v(z_0)$ for some $z_0\in\Omega$. Let
        $\psi\in\calE_0(\Omega)\cap C^{\infty}(\Omega)$ be a strictly
        plurisubharmonic exhaustion function 
        and let $s>0$ be such that $u(z_0)<s\psi(z_0) + v(z_0).$ Corollary 3.6 in
        \cite{Ceg00} gives, with $A= \{u(z)<s\psi(z) + v(z)\}$,
        $$\int_{A}(dd^c(s\psi + v))^n \leq \int_{A}(dd^cu)^n = 0.$$
        Hence $s^n\int_{A}(dd^c\psi)^n =0$ which implies that $A$ has
        Lebesgue measure $0$.
        Since the functions involved are plurisubharmonic, this means that
        $A=\emptyset$.
        This is a contradiction and the lemma is proved.
\end{proof}

\begin{proof}[Proof of Proposition \ref{prop_compare}]
        Let $K=\text{supp}\,(dd^c u)^n$. Since $K$ is compact,
        and since $u$ and $v$ are bounded upper semicontinuous functions,
        $\alpha>0$ may be chosen such that $\alpha v\le u$ on $K$.
        It then follows from Lemma \ref{lemma1}
        that $\alpha v\le u$ holds on all of $\Omega$.
        Similarly, there is $\beta>0$ such that $\beta u \le v$ on $\Omega$.
        Then Proposition \ref{ineq_prop} implies that
        $\mu_{\alpha^{-1}u}\le \mu_v\le\mu_{\beta u}$.
        Hence, if we let $a=\alpha^{-n}$ and $b=\beta^n$,
        the proposition follows.
\end{proof}

\begin{corollary}\label{cor_supp_S}
        There is a set $S\subset\dOmega$ such that
        $\text{supp}\,\mu_u=S$ for each $u\in\calF(\Omega)$, $u\ne 0$.
\end{corollary}

\begin{proof}
        Choose a function $v_0\in\calE_0(\Omega)$ with
        $\text{supp}\,(dd^c v_0)^n\subset\subset\Omega$,
        and let $S=\text{supp}\,\mu_{v_0}$.
        Let $u$ be an arbitrary function in $\calF(\Omega)$. Choose a sequence
        $\{u_j\}\subset\calE_0(\Omega)$ such that $u_j\searrow u$
        and $\text{supp}\,(dd^c u_j)^n\subset\subset\Omega$.
        Then Proposition \ref{prop_compare} implies that
        $\text{supp}\,\mu_{u_j}=S$ for each $j$. Moreover,
        $\mu_{u_1}\le\mu_{u_2}\le\cdots\le\mu_u$ and
        $\mu_{u_j}$ tends weak* to $\mu_u$, by Proposition \ref{ineq_prop}
        and Proposition \ref{appr_prop}.
        Hence $\text{supp}\,\mu_u= S$.
\end{proof}

Note that if $\mu$ is in the the weak* closure of
$\{\mu_u : u\in \calF(\Omega)\}$, then $\text{supp}\,\mu\subset S$.
Hence if $\Omega$ is B-regular, then the support set $S$ has to be all of
$\dOmega$, because of Theorem \ref{th3}.

On the other hand, if
$\Omega=\omega_1\times\omega_2\subset\mathbb{C}^n=\mathbb{C}^{n_1+n_2}$,
where $\omega_1\subset\mathbb{C}^{n_1}$ and
$\omega_2\subset\mathbb{C}^{n_2}$ are bounded hyperconvex domains,
then $S\subset\partial\omega_1\times\partial\omega_2$.
To see this, consider the function $u(z,w)=\max\,\{g_1(z),g_2(w)\}$ where
$g_k$ is the pluricomplex Green's function for $\omega_k$ with pole at some point
in $\omega_k$. Note that $g_k$ is continuous outside the pole and tends to
zero at the boundary of $\omega_k$.
Then $u\in\calF(\Omega)$ and
$\text{supp}\,(dd^c u)^n\subset\{(z,w)\in\Omega: g_1(z)=g_2(w)\}$.
Choose a sequence $\{\veps_j\}$ such that $\veps_j\searrow 0$.
Then $\Omega_j=\{(z,w)\in\Omega:u(z,w)<-\veps_j\}$ defines a fundamental
sequence of $\Omega$ and
$u^j:=\sup\,\{\varphi\in PSH(\Omega):
\varphi|_{\Omega\setminus\Omega_j}\le u|_{\Omega\setminus\Omega_j}\}=
\max\,\{u,-\veps_j\}$.
It follows that $\text{supp}\,(dd^c u^j)^n\subset
\{(z,w)\in\Omega: g_1(z)=g_2(z)\ge-\veps_j\}$, which implies that
$\text{supp}\,\mu_u\subset\partial\omega_1\times\partial\omega_2$.
Hence the claim follows from Corollary \ref{cor_supp_S}.

Using a similar argument, the following example shows that when
$\Omega=\mathbb{D}\times\mathbb{D}\subset\mathbb{C}^2$, then we have equality,
$S=\partial\mathbb{D}\times\partial\mathbb{D}$.

\begin{example}\label{ex_bidisc}
        Let $\Omega$ be the unit bidisc $\mathbb{D}\times\mathbb{D}$
        in $\mathbb{C}^2$. Then $\text{supp}\,\mu_u$ is equal to the
        distinguished boundary
        $\partial\mathbb{D}\times\partial\mathbb{D}$
        for each $u\in\calF(\Omega)$, $u\ne 0$.
        This follows from Corollary \ref{cor_supp_S}, if we
        for example consider the pluricomplex Green's function
        $g$ for $\Omega$ with pole at the origin.
        We then have that
        $g(z,w)=m\cdot\max\,\{\log|z|,\log|w|\}$, where the constant
        $m>0$ is chosen such that $\int_\Omega(dd^c g)^n = 1$.
        This is a function in $\calF(\Omega)$, and we can compute $\mu_g$
        explicitly.
        For $j=1,2,\ldots$, let $\Omega_j=\{(z,w):|z|<r_j, |w|<r_j\}$
        where $r_j=1-\frac{1}{j}$.
        Then $g^j:= \sup\,\{\varphi\in PSH(\Omega):
        \varphi|_{\Omega\setminus\Omega_j}\le
        g|_{\Omega\setminus\Omega_j}\}=
        m\cdot\max\,\{\log|z|,\log|w|, \log(r_j)\}$, from which it follows
        that $(dd^c g^j)^2 =
        m^2\cdot dd^c(\max\,\{\log|z|,\log(r_j)\})\wedge
        dd^c(\max\,\{\log|w|,\log(r_j)\})$.
        Since $\int_\Omega(dd^c g^j)^2=1$ for each $j$
        (see Section \ref{sec_constr}), we can conclude that
        $(dd^c g^j)^2=\sigma_j\times\sigma_j$, where $\sigma_j$
        is the normalized Lebesgue measure on the the circle
        $\partial \mathbb{D}(0,r_j)$.
        This implies that $\mu_g=\sigma\times\sigma$, where $\sigma$
        is the normalized Lesbegue measure on the unit circle.
\end{example}

\begin{remark}\label{rem_supp}
        Recall from Remark \ref{rem_pi} at the end of Section \ref{sec_constr}
        that Corollary \ref{cor_supp_S} and (\ref{pi_formula}) together give the
        partial integration formula
        \begin{equation}\label{pi_impl}
                h|_S=0\ \Rightarrow\ 
                \int_\Omega h\,(dd^c u)^n=
                \int_\Omega u\,dd^c h\wedge (dd^c u)^{n-1}.
        \end{equation}
        The implication (\ref{pi_impl})
        holds true for $h\in PSH(\Omega)\cap C(\bar\Omega)$ if
        $u\in\calF(\Omega)$,
        and for $h\in PSH(W)\cap L^\infty(W)$, $W\supset\bar\Omega$,
        if $u\in\calF^{a}(\Omega)$ (using Theorem \ref{th_equal} of Section
        \ref{sec_bdryval}).
        Here $S$ is the support set defined in Corollary \ref{cor_supp_S}.
        
        Furthermore, (\ref{J_ineq}) implies that
        \begin{equation}\label{S_ineq}
                \sup_\Omega\,h \le \sup_S\,h,
                \quad\forall\,h\in PSH^-(\Omega)\cap C(\bar\Omega).
        \end{equation}
        To see this, let $h\in PSH^-(\Omega)\cap C(\bar\Omega)$ be given.
        For $z\in\Omega$ fixed, let $g_z$ be the pluricomplex Green's
        function for $\Omega$ with pole at $z$.
        Then $(dd^c g_z)^n=\delta_z$ and we have
        that $h(z) = \int_\Omega h\,(dd^c g_z)^n\le
        \int_\dOmega h\,d\mu_{g_z}\le \sup_S\,h$.
        By the same argument, (\ref{S_ineq})
        holds true if $h$ is an upper bounded function in
        $PSH(W)$, where $W\supset\bar\Omega$.
\end{remark}

\begin{remark}\label{rem_Henkin}
        Another property of the measures $\mu_u$ is that they are so called
        \textit{Henkin measures}
        (a kind of measure introduced by Henkin in \cite{Hen}).
        This means that $$\lim_{k\to\infty}\int_\dOmega f_k\,d\mu_u=0$$ for each uniformly
        bounded sequence $\{f_k\}$ in $A(\Omega)$ such that $\lim_{k\to\infty} f_k(z)=0$
        for all $z\in\Omega$. Here $A(\Omega)$ denotes the functions
        that are holomorphic on $\Omega$ and continuous on $\bar\Omega$.
        To see that this holds, take such a sequence $\{f_k\}$ and let
        $\{\varphi_k\}=\{\text{Re}\,f_k\}$. From (\ref{limit_1}) and
        (\ref{eq_1}) it follows that
        $$\lim_{k\to\infty}\int_\dOmega\varphi_k\,d\mu_u =
        \lim_{k\to\infty}\left(\lim_{j\to\infty} \int_\Omega\varphi_k\,(dd^c u^j)^n\right) =
        \lim_{k\to\infty}\int_\Omega\varphi_k\,(dd^c u)^n = 0$$
        for each $u\in\calF(\Omega)$,
        since $\varphi_k$ is uniformly bounded and
        $\int_\Omega (dd^c u)^n<\infty$.
        Since the same holds for $\{\psi_k\}=\{\text{Im}\,f_k\}$, it follows
        that $\lim_{k\to\infty}\int_\dOmega f_k\,d\mu_u=0$.
               
        This property can be used to show the following fact
        about the support of the measures $\mu_u$.
        Suppose that $u\in\calF(\Omega)$ and that $K\subset\dOmega$ is a peak
        set for $A(\Omega)$.
        Let $f\in A(\Omega)$ be a peak function for $K$ and define
        $f_k(z) = (f(z))^k$, for $z\in\bar\Omega$ and $k=1,2,\ldots$.
        Then $\{f_k\}$ satisfies the assumptions above, so
        $\lim_{k\to\infty} \int_\dOmega f_k\,d\mu_u = 0$.
        But we also have that
        $\lim_{k\to\infty} \int_\dOmega f_k\,d\mu_u = \mu_u(K)$. Hence
        $\mu_u(K)=0$ for each peak set $K$ and each $u\in\calF(\Omega)$.
\end{remark}


\bigskip\bigskip
\section{Boundary values}\label{sec_bdryval}

\noindent
In this section we define and study boundary values of plurisubharmonic
funtions, with respect to the measures $\mu_u$.

\begin{lemma}\label{bdry_lemma}
        Assume that $u\in\calF(\Omega)$ and
        $g\in PSH(\Omega)\cap L^\infty(\Omega)$.
        Then $\{g\,(dd^c u^j)^n\}$ is weak*-convergent.
\end{lemma}

\begin{proof}
         By the same argument as in Theorem \ref{th1} it is enough to prove
         that the limit $\lim_{j\to\infty}
         \int_\Omega \varphi g\,(dd^c u^j)^n$ exists
         for all $\varphi\in PSH^-(\Omega)\cap L^\infty(\Omega)$.
         Given such a function $\varphi$, take $M,N\ge 0$ such that
         $\varphi+M\ge0$ and $g+N\ge 0$.
         Then $(\varphi+M)^2$, $(g+N)^2$, $(\varphi+M+g+N)^2\in
         PSH(\Omega)\cap L^\infty(\Omega)$, so if $\psi$ is any of these then
         $\lim_{j\to\infty} \int_\Omega \psi\,(dd^c u^j)^n$ exists by Theorem \ref{th1}.
         Expanding $((\varphi+M)+(g+N))^2$, it follows that the limit exists
         for $\psi=(\varphi+M)(g+N)$ and then finally for $\psi=\varphi g$
         (using Theorem \ref{th1} again).
\end{proof}

Using this lemma, together with standard measure theory, we can make the
following definition.

\begin{definition}\label{b_val_def}
        For $u\in\calF(\Omega)$ and $g\in PSH(\Omega)\cap L^\infty(\Omega)$,
        let $g^u$ be the function in $L^\infty(\dOmega,\mu_u)$ such that
        $\lim_{j\to\infty} g\,(dd^c u^j)^n = g^u\,d\mu_u$.
\end{definition}

We may consider $g^u$ as the boundary values of $g$ with respect to $\mu_u$.
Note that, at least formally, $g^u$ depends on both $g$ and $u$.
However, the following theorems describe some situations when this definition
agrees with other notions of boundary values.

\begin{theorem}\label{th_equal}
        Assume that $u\in\calF^{a}(\Omega)$ and 
        $g\in PSH(W)\cap L^\infty(W)$ where $W$ is a
        bounded domain containing $\bar\Omega$.
        Then $g^u=g|_\dOmega$ a.e.\ $(\mu_u)$.
\end{theorem}

\begin{proof}
        Note that if $M$ is a constant then $(g-M)^{u}=g^{u}-M$,
        so we may assume that $g\le 0$.
        Let $t\in C(\bar\Omega)$, $t\ge 0$ be given. Then it follows, in the
        same way as in the proof of Lemma \ref{tec_lemma}, that
        $$\int_\dOmega tg^u\,d\mu_u =
        \lim_{j\to\infty}\int_\Omega tg\,(dd^c u^j)^n \le
        \int_\dOmega tg\,d\mu_u.$$
        Thus $g^u\le g$ a.e.\ $(\mu_u)$, so it remains to prove
        that $\int_\dOmega g^u\,d\mu_u = \int_\dOmega g\,d\mu_u$.
        Choose $K$ such that
        $\Omega\subset\subset K \subset\subset W$.
        Given $\veps > 0$ there is an open set $U_\veps\subset W$ and a function
        $g_\veps\in C_0(W)$ such that $\inf_W g\le g_\veps\le 0$,
        the relative capacity $cap\,(U_\veps, W) < \veps$ and
        $K\setminus U_\veps\subset\{z\in W: g(z) = g_\veps(z)\}$
        (for definition and properties of relative capacity, see \cite{B&T}).
        It follows that
        \begin{equation*}
        \begin{split}
                \int_\dOmega g^u\,d\mu_u & =
                \lim_{j\to\infty}\int_\Omega g\,(dd^c u^j)^n =
                \\
                & = \lim_{j\to\infty}\int_{\Omega\cap U_\veps} g\,(dd^c u^j)^n +
                \lim_{j\to\infty}\int_{\Omega\setminus U_\veps}
                g_\veps\,(dd^c u^j)^n\ge\\
                & \ge \lim_{j\to\infty}\int_{\Omega\cap U_\veps} g\,(dd^c u^j)^n +
                \int_\dOmega g_\veps\,d\mu_u =\\
                & = \lim_{j\to\infty}\int_{\Omega\cap U_\veps} g\,(dd^c u^j)^n +
                \int_{\dOmega\cap U_\veps} g_\veps\,d\mu_u + 
                \int_{\dOmega\setminus U_\veps} g\,d\mu_u\ge\\
                & \ge \lim_{j\to\infty}\int_{\Omega\cap U_\veps} g\,(dd^c u^j)^n +
                \int_{\dOmega\cap U_\veps} g_\veps\,d\mu_u + 
                \int_{\dOmega} g\,d\mu_u.
        \end{split}
        \end{equation*}
        Let $h_\veps=\sup\,\{\psi\in PSH^-(W): \psi|_{U_\veps}\le-1\}$,
        we then have that
        \begin{equation*}
        \begin{split}
                0&\ge \int_\dOmega g^u\,d\mu_u - \int_{\dOmega} g\,d\mu_u\ge\\
                & \le
                \lim_{j\to\infty}\int_{\Omega\cap U_\veps} g\,(dd^c u^j)^n +
                \int_{\dOmega\cap U_\veps} g_\veps\,d\mu_u\ge\\
                & \ge\left(\inf_W g\right)
                \left(\lim_{j\to\infty}\int_{\Omega\cap U_\veps}(dd^c u^j)^n
                + \int_{\dOmega\cap U_\veps}d\mu_u\right)=\\
                & =\left(-\inf_W g\right)
                \left(\lim_{j\to\infty}\int_{\Omega\cap U_\veps}h_\veps\,(dd^c u^j)^n
                + \int_{\dOmega\cap U_\veps}h_\veps\,d\mu_u\right)\ge\\
                & \ge\left(-\inf_W g\right)
                \left(\lim_{j\to\infty}\int_{\Omega}h_\veps\,(dd^c u^j)^n
                + \int_{\dOmega}h_\veps\,d\mu_u\right)\ge\\
                &\ge 2\left(-\inf_W g\right)
                \int_\Omega h_\veps\,(dd^c u)^n,
         \end{split}
         \end{equation*}
         where we have used (\ref{ineq_1}) and Lemma \ref{tec_lemma} in the
         last inequality. From Lemma 1.9 in \cite{CegKol}, using that
         $u\in\calF^{a}(\Omega)$ and that $cap\,(U_\veps, W) < \veps$,
         it follows that this last integral tends to zero as $\veps\searrow 0$,
         which completes the proof.
\end{proof}

The following theorem may be compared with the definitions in Section
\ref{sec_prel}.

\begin{theorem}\label{th2}
        Suppose that $H\in \calM(\Omega)\cap L^\infty(\Omega)$.
        Then, for every $u\in\calF^{a}(\Omega)$ and every $g\in
        \calF(\Omega,H)$ 
        such that $\int_\Omega g\,(dd^cu)^n > -\infty$,
        $g\,(dd^cu^j)^n$ is weak*-convergent to $H^{u}\,d\mu_u$.
\end{theorem}

\begin{proof}
        By the same argument as in Theorem
        \ref{th1}, it is enough to prove that
        $$
        \lim_{j\to\infty} \int_\Omega tg\,(dd^c u^j)^n = 
        \lim_{j\to\infty} \int_\Omega t H (dd^c u^j)^n,
        \quad\forall\,t\in PSH^-(\Omega)\cap L^\infty(\Omega).
        $$
        Since $g\in \calF(\Omega,H)$ there is a
        $\psi\in\calF(\Omega)$ such that
        $\psi + H \le g\le H$.
        We may assume that $\psi\ge g$ (otherwise, look at
        $\psi_0=\max\,\{\psi,g\}$). We may also (after dividing by suitable
        constants) assume that $-1\le t\le 0$ and $-1\le H\le 0$.
        Now,
        $$
        \int_\Omega tg\,(dd^c u^j)^n =
        \int_\Omega t(g-H)\,(dd^c u^j)^n +
        \int_\Omega t H (dd^c u^j)^n
        $$
        where
        $0\le \int_\Omega t(g-H)\,(dd^c u^j)^n = 
        \int_\Omega (-t)(H-g)\,(dd^c u^j)^n \le
        \int_\Omega (-t)(-\psi)\,(dd^c u^j)^n \le
        \int_\Omega (-\psi)\,(dd^c u^j)^n$.
        Using partial integration in $\calF(\Omega)$ we have the following
        \begin{equation*}
        \begin{split}
                \int_\Omega (-\psi)\,(dd^c u^j)^n & =
                \int_\Omega (-u^j)\,dd^c\psi\wedge(dd^c u^j)^{n-1} \le
                \int_\Omega (-u)\,dd^c\psi\wedge(dd^c u^j)^{n-1}=\\
                 & = \int_\Omega (-u^j)\,dd^c\psi\wedge dd^c u
                 \wedge(dd^c u^j)^{n-2} \le \ldots\le\\
                 &\le\int_\Omega (-u^j)\,dd^c\psi\wedge(dd^c u)^{n-1} =
                 I_j\le
                 \int_\Omega (-u)\,dd^c\psi\wedge(dd^c u)^{n-1}=\\
                 & =\int_\Omega (-\psi)\,(dd^c u)^n \le
                  \int_\Omega (-g)\,(dd^c u)^n <+\infty.
        \end{split}
        \end{equation*}
        Since $u^j$ increases to zero outside a pluripolar set and
        $dd^c\psi\wedge(dd^c u)^{n-1}$ vanishes on pluripolar sets
        (see Section \ref{sec_prel}, Lemma \ref{l_pp}), it follows
        that $I_j\searrow 0$ when $j\to+\infty$.
        This proves the theorem.
\end{proof}

\begin{remark}\label{rem_th2}
        If $g\in L^\infty(\Omega)$ then $\int_\Omega g\,(dd^cu)^n > -\infty$
        for every $u\in\calF(\Omega)$. Furthermore, $\psi\ge g$ implies that
        $\psi$ is bounded as well, so $dd^c\psi\wedge(dd^c u)^{n-1}$ vanishes
        on pluripolar sets for every $u\in\calF(\Omega)$ (Lemma \ref{l_pp}).
        Thus for bounded functions $g$ in $\calF(\Omega,H)$,
        the conclusion
        $g^{u}\,d\mu_u = H^{u}\,d\mu_u$
        holds for every $u\in\calF(\Omega)$.
\end{remark}

Suppose that we have a bounded plurisubharmonic function on $\Omega$ and want to
approximate it with
plurisubharmonic functions that are continuous on $\bar\Omega$.
The following theorem gives a condition for when this implies
weak*-convergence on the boundary.

\begin{theorem}
        Assume that
        $u\in \calF(\Omega)$ and  $\mu_u = \lim_{j\to\infty}(dd^cu^{j})^{n}$.
        Let $\{\varphi_j\}$ be a sequence in $PSH(\Omega)\cap
        C(\bar\Omega)$ such that $0\leq \varphi_j\leq 1$.    
        If $\varphi_j$
        tends to $\varphi\in PSH(\Omega)\cap L^\infty(\Omega)$
        in the sense of distributions, then
        $\varphi_j\,d\mu_u $ tends weak* to $\varphi^{u}\, d\mu_u$ if and only if
        $\lim_{j\to\infty}\int \varphi_j\,d\mu_u = \int\varphi^{u}\,d\mu_u$.
\end{theorem}

\begin{proof}
        By Corollary \ref{cor_max} we may assume that
        $u\in\calF^{a}(\Omega)$.
        The condition in the theorem is obviously necessary, we prove it is also
        sufficient.
        First, note that for $\{\psi_k\}\subset PSH(\Omega)\cap
        C(\bar\Omega)$, $\psi_k\ge 0$, the following holds.
        For $k$ fixed,
        $(\sup_{l\ge k}\psi_l)^*\in PSH(\Omega)\cap L^\infty(\Omega)$,
        therefore $(\sup_{l\ge k}\psi_l)^*\,(dd^c u^j)^n$ is weak*-convergent
        (as $j\to\infty$) by Lemma \ref{bdry_lemma}.
        Furthermore, since $(\sup_{l\ge k}\psi_l)=(\sup_{l\ge k}\psi_l)^*$
        outside a pluripolar set and $u^j\in\calF^{a}(\Omega)$
        (since $u\in\calF^{a}(\Omega)$), the star may be removed.
        We claim that
        \begin{equation}\label{part_i}
                \lim_{j\to\infty}\: (\sup_{l\ge k}\psi_l)\,(dd^c u^j)^n=
                (\sup_{l\ge k}\psi_l)\,d\mu_u.
        \end{equation}
        Given $f\in C(\bar\Omega)$, $f\ge0$ it follows from (\ref{limit_1})
        that for each $m$
        $$
        \lim_{j\to\infty}\int_\Omega f(\sup_{l\ge k}\psi_l)\,(dd^c u^j)^n \ge
        \lim_{j\to\infty}\int_\Omega f(\sup_{m\ge l\ge k}\psi_l)\,(dd^c u^j)^n=
        \int_\dOmega f(\sup_{m\ge l\ge k}\psi_l)\,d\mu_u,
        $$
        where the last integral tends to
        $\int_\dOmega f(\sup_{l\ge k}\psi_l)\,d\mu_u$ as $m\to\infty$.
        It follows that
        $\lim_{j\to\infty}\: (\sup_{l\ge k}\psi_l)\,(dd^c u^j)^n\ge
        (\sup_{l\ge k}\psi_l)\,d\mu_u$.
        On the other hand, by (\ref{ineq_1}) and (\ref{limit_1})
        $$\int_\Omega (\sup_{m\ge l\ge k}\psi_l)\,(dd^c u^j)^n \le
        \int_\dOmega (\sup_{m\ge l\ge k}\psi_l)\,d\mu_u$$
        for each $m$ and $j$. So by letting $m\to\infty$ we have that
        $\int_\Omega (\sup_{l\ge k}\psi_l)\,(dd^c u^j)^n \le
        \int_\dOmega (\sup_{ l\ge k}\psi_l)\,d\mu_u$,
        which proves the claim.
        
        Now, let $\{\varphi_{j_m}\,d\mu_u\}$ be any weak*-convergent subsequence
        of $\{\varphi_j\,d\mu_u\}$. (Such a sequence exists by the same
        reasoning as in the proof of Theorem \ref{th3}.)
        Then, by standard measure theory,
        the limit measure is equal to $\varphi_0\,d\mu_u$
        for some $\varphi_0\in L^\infty(\mu)$.
        We will show that $\varphi_0=\varphi^u$ a.e.\ $(\mu)$.
        It then follows that the original sequence itself converges to
        $\varphi^{u}\,d\mu_u$, and the proof will be complete.
        
        From $L^2$-theory it follows that we may choose
        $\psi_k=\frac{1}{M_k}\sum_{l=1}^{M_k}\varphi_{j_{m_l}}$
        such that $\psi_k\to\varphi_0$ in $L^2(\mu)$ and then a subsequence
        converging to $\varphi_0$ a.e.\ $(\mu)$, for simplicity call it
        $\{\psi_k\}$.
        Since by assumption the original sequence $\{\varphi_j\}$ tends
        to $\varphi$ in the sense of distributions, the same holds for
        $\{\psi_k\}$.
        Now, for $f\in C(\bar\Omega)$, $f\ge 0$, using the definition of
        $\varphi^u$, (\ref{part_i}) and monotone convergence,
        \begin{equation*}
        \begin{split}
                \int_\dOmega f\varphi^u\,d\mu_u & = 
                \lim_{j\to\infty}\int_\Omega f\varphi\,(dd^c u^j)^n =
                \text{(Lemma 1.4 in \cite{CegKol})} =\\
                & = \lim_{j\to\infty}\lim_{k\to\infty} \int_\Omega f\psi_k\,(dd^c u^j)^n \le
                \lim_{k\to\infty}\lim_{j\to\infty} \int_\Omega f(\sup_{l\ge k}\psi_l)\,
                (dd^c u^j)^n =\\
                & = \lim_{k\to\infty}\int_\dOmega f(\sup_{l\ge k}\psi_l)\,d\mu_u = 
                \int_\dOmega f (\limsup_{k\to\infty}\psi_k)\,d\mu_u.
        \end{split}
        \end{equation*}
        From this it follows that $\varphi^u\le\limsup_{k\to\infty}\psi_k$
        a.e.\ $(\mu)$, which implies that $\varphi^u\le\varphi_0$ a.e.\ $(\mu)$.
        Furthermore, $\int_\dOmega\varphi_0\,d\mu_u =
        \lim_{m\to\infty} \int_\dOmega\varphi_{j_m}\,d\mu_u =
        \int_\dOmega\varphi^u\,d\mu_u$, by assumption, so
        $\varphi^u=\varphi_0$ a.e.\ $(\mu)$.
        Hence the theorem is proved.
\end{proof}


\bigskip\bigskip
\section{More boundary measures}\label{sec_more}

\noindent
Let $\nu$ be a positive measure on $\Omega$ with finite total mass. Then there
is a positive measure $\mu\ne 0$ which is supported by $\dOmega$, vanishes on
pluripolar sets and such that%
\begin{equation}\label{mes_ineq}
        \int_\Omega\varphi\,d\nu \leq \int_\dOmega\varphi\,d\mu,
        \quad\forall\,\varphi\in PSH^-(\bar\Omega),
\end{equation}
where $PSH^-(\bar\Omega)
=\{\varphi: \varphi\in PSH^-(\Omega'),\, \Omega'\supset\bar\Omega\}$.
To see this, let $P_\nu$ denote the pluricomplex potential of $\nu$ relative
to $\Omega$, i.e.\ $P_\nu(z)=\int_\Omega g(z,w)\,d\nu(w)$,
where $g(z,w)$ is the pluricomplex Green's function for $\Omega$ with
pole at $w$.
Then Theorem 1.1 in \cite{Ceg4} says that $P_\nu\in\calF(\Omega)$ and that
$$\int_\Omega-\varphi\,(dd^c P_\nu)^n\le
\left(\nu(\Omega)\right)^{n-1}\int_\Omega -\varphi\,d\nu,
\quad\forall\,\varphi\in PSH^-(\Omega).$$
Moreover,
$\int_\Omega\varphi\,(dd^c P_\nu)^n\le
\int_\dOmega\varphi\,d\mu_{P_\nu}$ for each $\varphi\in PSH^-(\bar\Omega)$,
by Remark \ref{rem_pi} at the end of Section \ref{sec_constr}.
Hence, the claim follows if we take $\mu = \left(\nu(\Omega)\right)^{-n+1}\mu_{P_\nu}$.

Conversely, if a positive measure $\mu$ on $\dOmega$ is such that
(\ref{mes_ineq}) holds for some finite measure $\nu$ on $\Omega$,
we would like to find an approximation procedure, similar to the one in Section
\ref{sec_constr}. A motivation is that we are interested in boundary values of
plurisubharmonic functions with respect to $\mu.$

We will study the case when $\nu$ vanishes on all pluripolar subsets of $\Omega$
and $\Omega$ belongs to a more restricted class of hyperconvex domains: 
\begin{enumerate}[(\ref{sec_more}a)]
        \item\label{cond_a}
                $\Omega$ and $\{\Omega_k\}$
                are hyperconvex domains with
                $\Omega\subset\subset\Omega_{k+1}\subset\subset\Omega_k$,
                such that for each $t\in\calF(\Omega)$ there is a sequence
                $\{t_k\}$, where $t_k\in\calF(\Omega_k)$ and
                $t_k\nearrow t$ a.e.\ on $\Omega$.
        \item\label{cond_b}
                $\Omega$ is not thin at any of its boundary points, so
                that $\limsup_{\Omega\ni z\to\xi}v(z)=v(\xi)$ for
                each $\xi\in\dOmega$ if
                $v\in PSH^-(\bar\Omega)$.
\end{enumerate}
Conditions for the approximation property in (\ref{sec_more}\ref{cond_a}) to hold
true have been studied in for example \cite{Benel} and \cite{CegHed}.
Examples of domains satisfying (\ref{sec_more}\ref{cond_a}) and 
(\ref{sec_more}\ref{cond_b}) are polydiscs and strictly pseudoconvex domains.
Note that if $t$ is bounded, we may assume that each $t_k$ is bounded.
 
 \begin{theorem}\label{th_more}
        Let $\Omega$ be a domain satisfying
        (\ref{sec_more}\ref{cond_a}) and (\ref{sec_more}\ref{cond_b}).
        Assume that $\mu$ is a positive measure on $\partial\Omega$, vanishing
        on pluripolar sets.
        Then there is a sequence $\{w_k\}$ in $\calF^{a}(\bar\Omega) = 
        \{u: u\in \calF^{a}(\Omega'),\, \Omega'\supset\bar\Omega\}$
        such that
        $\text{supp}\,(dd^c w_k)^n\subset\subset\Omega$,
        $\int_\Omega (dd^c w_k)^n\le \int_\dOmega d\mu$, and
        $(dd^cw_k)^n$ tends weak* to $\mu$ as $ k\to \infty$.
        
        Furthermore, if there is a finite positive measure $\nu$ on
        $\Omega$, vanishing on pluripolar sets, such that (\ref{mes_ineq})
        holds, then
        $\lim_{k\to\infty}\int_\Omega t\,(dd^cw_k)^n = 0$ for each
        $t\in\calF(\Omega)\cap L^\infty(\Omega)$.
        Hence $t\,(dd^cw_k)^n$ tends weak* to $0$
        for each $t\in\calF(\Omega)\cap L^\infty(\Omega)$.
\end{theorem}
 
If we compare this theorem with the results in the previous sections we have the
following. In the setting of Section \ref{sec_constr} we know that if
$u\in\calF^{a}(\Omega)$ and $\varphi \in PSH(W)\cap L^\infty(W)$,
$W\supset\bar\Omega$, then
$\int_\Omega\varphi\,(dd^cu)^n\leq \int_\Omega\varphi\,(dd^cu^j)^n$
which increases to $\int_\dOmega\varphi\, d\mu_u$ (see Theorem
\ref{th_equal}).
In particular it follows that when $\mu=\mu_u$ for some $u\in\calF^{a}(\Omega)$,
then (\ref{mes_ineq}) is satisfied if we take $\nu=(dd^c u)^n$.
We also have that $\int_\Omega (dd^c u^j)^n=\int_\dOmega d\mu_u$ and
$\lim_{j\to\infty} t\,(dd^c u^j)^n=0$ for each $t\in\calF(\Omega)\cap L^\infty(\Omega)$
(see Remark \ref{rem_th2}).
Hence, the approximation procedure in Theorem \ref{th_more} is similar to the one
in the previous sections, and it applies to a larger class of boundary
measures, see also Example \ref{ex_more}.

\begin{lemma}\label{lemma2}
        Let ${\{\mu^{j}_k\}}_{j,k}$ be a sequence of positive measures on
        $\bar\Omega$ with uniformly bounded mass.
        Suppose that, for each fixed $k$, $\mu^{j}_k$ tends weak* to $\mu$ as
        $j\to\infty$.
        Then there is a subsequence ${\{\mu^{j_k}_k\}}_k$ such that
        $\mu^{j_k}_k$ tends weak* to $\mu$ as $k\to \infty$.
\end{lemma}

\begin{proof}
        Let $\{t_l\}$ be a dense sequence in $C(\bar\Omega)$. For each $k$ we
        choose $j_k$ such that
        $$\left|\int_{\bar\Omega} t_l\,d\mu -
        \int_{\bar\Omega} t_l\,d\mu^{j_k}_k\right| <
        \frac{1}{k},\quad 1\leq l \leq k.$$
        It follows that $\mu^{j_k}_k$ tends weak* to $\mu$ as $k\to \infty$,
        since $\{t_l\}$ is dense and the measures have uniformly bounded
        total mass.
\end{proof}

\begin{proof}[Proof of Theorem \ref{th_more}]        
        For each $k$, the measure $\mu$ can be regarded as a finite
        measure on
        $\Omega_k$ which vanishes on pluripolar sets. Hence there is
        $u_k\in\calF^{a}(\Omega_k)$ such that $(dd^cu_k)^n = \mu$,
        see Lemma 5.14 in \cite{Ceg2}.
        Choose a fundamental sequence
        $\{\omega_j\}$ of $\Omega$, i.e.\  
        $\omega_j\subset\subset \omega_{j+1}\subset\subset\Omega$
        and $\cup_{j=1}^\infty\omega_j=\Omega$.
        For each $k$ and $j$, define
        $u_k^j = \sup\,\{\varphi\in PSH^-(\Omega_k):
        \varphi|_{\omega_j}\leq u_k|_{\omega_j}\}$.
        Then $u_k^j\in\calF^{a}(\Omega_k)$
        (note that $(u_k^j)^*=u_k^j$ since $\omega_j$ is open,
        so $u_k^j$ is plurisubharmonic) and we have the
        following:
        \begin{enumerate}[(i)]
                \item\label{stat_1}
                $\text{supp}\,(dd^cu_k^j)^n \subset \partial\omega_j$,
                $u_k^j \geq u_k$ on $\Omega_k$,
                $\int_{\Omega_k}(dd^cu_k^j)^n \leq 
                \int_{\Omega_k}(dd^cu_k)^n =\int_\dOmega d\mu$.
                \item\label{stat_2}
                If $j_1\leq j_2$ then $u_k^{j_1} \geq u_k^{j_2}$ on $\Omega_k$.
                \item \label{stat_3}
                $\lim_{j \to \infty}u_k^j = u_k$ on $\Omega_k$.
        \end{enumerate}
        The first two statements are obvious.
        For the proof of the third, let $v_k=\lim_{j}u_k^j$.
        Then $v_k\in\calF(\Omega_k)$, $v_k\ge u_k$ on $\Omega_k$ and
        $v_k=u_k$ on $\Omega$. Thus $v_k(\xi)=u_k(\xi)$ for $\xi\in\dOmega$,
        using the assumption (\ref{sec_more}\ref{cond_b}),
        so $v_k\le u_k$ on $\Omega_k$
        by Lemma \ref{lemma1} and the statement follows.
        Now, (\ref{stat_2}) and (\ref{stat_3}) imply that
        $(dd^c u_k^j)^n$ tends weak* to $(dd^c u_k)^n=\mu$ as $j\to\infty$,
        for each fixed $k$.
        Hence, by (\ref{stat_1}) we can
        use Lemma \ref{lemma2} to pick $\{j_k\}$ such that
        $(dd^cu_k^{j_k})^n$
        tends weak* to $\mu$ as $k\to\infty$. This completes the first part of
        the theorem, if we let $w_k = u_k^{j_k}$.
                
        It remains to prove that $\lim_{k\to \infty}\int
        t\,(dd^cu_k^{j_k})^n = 0$ for all
        $t\in\calF(\Omega)\cap L^\infty(\Omega)$, assuming that
        (\ref{mes_ineq}) holds.
        Given $t\in\calF(\Omega)\cap L^\infty(\Omega)$
        there is by (\ref{sec_more}\ref{cond_a})
        a sequence
        $\{t_k\}$ with $t_k\in\calF(\Omega_k)\cap L^\infty(\Omega_k)$
        such that
        $t_k$ increases a.e.\ to $t$ on $\Omega$.
        Now, 
        $$\int_{\Omega_k} t\,(dd^cu_k^{j_k})^n\geq
        \int_{\Omega_k} t_k\,(dd^cu_k^{j_k})^n \geq
        \int_{\Omega_k} t_k\,(dd^cu_k)^n =
        \int_\dOmega t_k\, d\mu \geq
        \int_\Omega t_k\,d\nu>-\infty$$
        so it follows that
        $$\liminf_{k\to \infty}
        \int_{\Omega_k} t\,(dd^cu_k^{j_k})^n
        \geq \int_\Omega t\,d\nu.$$
        Define $t^{i}= \sup\,\{\varphi\in PSH(\Omega):
        \varphi|_{\Omega\setminus\omega_i}\le
        t|_{\Omega\setminus\omega_i}\}$.
        Then $t^{i}\in\calF(\Omega)\cap L^\infty(\Omega)$ and
        $t^{i}=t$ on $\Omega\setminus\omega_i$, so
        $$\liminf_{k\to \infty}\int_{\Omega_k} t\,(dd^cu_k^{j_k})^n =
        \liminf_{k\to \infty}\int_{\Omega_k} t^{i}\,(dd^cu_k^{j_k})^n\geq
        \int_\Omega t^{i}\,d\nu,$$
        by the above calculations.
        Now, the left hand side is independent of $i$, while the right hand
        side tends to $0$ when $i$ tends to $\infty$,
        since $\nu$ vanishes on pluripolar sets.
        This completes the proof.
\end{proof}

The reason not to keep $k$ fixed in the proof above, is to be able to prove
the second part of the theorem. Also, one can prove that
$\lim_{k\to\infty} u_k^j=0$ a.e.\ on $\Omega$, for each fixed $j$.

\begin{remark}
        Suppose that $v\in PSH^-(\Omega)$ satisfies
        $\tilde v \ge v \ge \tilde v + \psi$ for some
        $\psi\in\calF(\Omega)\cap L^\infty(\Omega)$ and
        that $\tilde v\in C(\bar\Omega)$.
        (Thus, $v$ is a function in $\calF(\Omega, \tilde v)$
        with some additional properties, see Section \ref{sec_prel}.)
        Then the preceeding theorem implies that
        \begin{equation}\label{limit_6}
                \lim_{k\to\infty} v\,(dd^c u_k^{j_k})^n = \tilde v\,d\mu,
        \end{equation}
        where the limit is in weak* sense.
        To see this, take $f\in C(\bar\Omega)$,
        $f\ge 0$. Then by the theorem we have that
        $\lim_{k\to\infty}\int_\Omega f\tilde v\,(dd^c u_k^{j_k})^n =
        \int_\dOmega f\tilde v\,d\mu$ and that
        $0\ge\int_\Omega f\psi\,(dd^c u_k^{j_k})^n\ge
        {\max f}\cdot\int_\Omega\psi\,(dd^c u_k^{j_k})^n$,
        where the last integral tends to $0$ as $k\to\infty$.
        Hence the inequality $f\tilde v \ge fv \ge f\tilde v + f\psi$
        implies that $\lim_{k\to\infty}\int_\Omega fv\,(dd^c u_k^{j_k})^n =
        \int_\dOmega f\tilde v\,d\mu$, and (\ref{limit_6}) follows.
        
        Furthermore, if we assume 
        that $\int_\Omega \varphi\,d\nu > -\infty$
        for all $\varphi\in\calF(\bar\Omega)$,
        then (\ref{limit_6}) holds for all
        $v\in\calF(\Omega,\tilde v)$ where $\tilde v\in C(\bar\Omega)$.
        This is due to the fact that the boundedness of $t$ in the second
        part of Theorem \ref{th_more} is used only to ensure that
        $\int_\Omega t_k\,d\nu>-\infty$ (because if $t$ is bounded then
        $t_k$ is bounded). Hence the assumption that
        $t\in\calF(\Omega)\cap L^\infty(\Omega)$ can be replaced by
        the assumption that
        $t\in\calF(\Omega)$ and $\int_\Omega \varphi\,d\nu > -\infty$
        for all $\varphi\in\calF(\bar\Omega)$.
\end{remark}

\begin{example}\label{ex_more}
        Let $\Omega$ be the unit bidisc $\mathbb{D}\times\mathbb{D}$
        in $\mathbb{C}^2$. Let $\mu$ and $\nu$ be defined by
        \begin{equation*}
                \mu=\sigma_1\times dV_\frac{1}{2}\ \ \text{and}\ \ 
                \nu=\sigma_\frac{1}{2}\times dV_\frac{1}{2},
        \end{equation*}
        where $\sigma_r$ denotes the normalized Lebesgue measure on
        the circle $\partial\mathbb{D}(0,r)$ and $dV_\frac{1}{2}$ the
        normalized Lebesgue measure on the
        disc $\mathbb{D}(0,\frac{1}{2})$.
        Then $\mu$ and $\nu$ satisfies (\ref{mes_ineq}), so Theorem
        \ref{th_more} tells us that we can approximate $\mu$ from
        the inside of $\Omega$ by our procedure.
        Moreover, by Example \ref{ex_bidisc} we see that $\mu$
        is not in the weak* closure of $\{\mu_u: u\in\calF(\Omega)\}$.
        Hence, we do reach more measures by the method in this
        section than we could before.
\end{example}


\end{document}